\pdfoutput=1
\documentclass[12pt]{amsart}

\usepackage[top=22mm, bottom=22mm, left=30mm, right=30mm]{geometry}
\usepackage[export]{adjustbox}
\usepackage{subfig}

\usepackage{caption}
\usepackage{wrapfig}
\usepackage{pinlabel, color}
\usepackage{multirow}

\usepackage{amsthm} %%% Defines \theoremstyle
\usepackage{amsmath} %%% Defines \numberwithin, \operatorname
\usepackage{amssymb} %%% Defines \mathbb

\usepackage{microtype} 
\usepackage{pinlabel} 

\usepackage[
hidelinks,
pagebackref]{hyperref}

\renewcommand*{\backref}[1]{}
\renewcommand*{\backrefalt}[4]{
  \ifcase #1 
  [No citations.]
  \or [#2]
  \else [#2]
  \fi }

\let\originalleft\left
\let\originalright\right
\renewcommand{\left}{\mathopen{}\mathclose\bgroup\originalleft}
\renewcommand{\right}{\aftergroup\egroup\originalright}

\numberwithin{equation}{section} 

\theoremstyle{plain}

\theoremstyle{definition}

\newtheorem{exercise}[equation]{Exercise}

\newtheorem*{question*}{Question}
\newtheorem*{answer*}{Answer}
\newtheorem*{application*}{Application}

\theoremstyle{remark}

\newtheorem*{remark*}{Remark}

\newtheorem*{case*}{Case}
\newtheorem*{step*}{Step}

\newtheorem*{claim*}{Claim}

\newcommand{\reffig}[1]{Figure~\ref{Fig:#1}}

\newcommand{\fakeenv}{} %%% prints the emptystring

{ 
 \renewcommand{\fakeenv}{#2} %%% So now \fakeenv prints #2
 \theoremstyle{plain} 
 \newtheorem*{\fakeenv}{#1~\ref{#2}} %%% so now #2 is the name of a
 \begin{\fakeenv}
}
{
 \end{\fakeenv}
}

\newenvironment{restated}[2]  
{ 
 \renewcommand{\fakeenv}{#2} 
 \theoremstyle{definition} 
 \newtheorem*{\fakeenv}{#1~\ref{#2}} %%% so now #2 is the name of a
 \begin{\fakeenv}
}
{
 \end{\fakeenv}
}

\begin{document}

\title[Conformally correct]{Conformally correct tilings}

\author[Schleimer]{Saul Schleimer}
\address{\hskip-\parindent
        Mathematics Institute\\
        University of Warwick\\
        Coventry CV4 7AL, United Kingdom}
\email{s.schleimer@warwick.ac.uk}
%\urladdr{http://www.math.rutgers.edu/$\sim$saulsch}

\author[Segerman]{Henry Segerman}
\address{\hskip-\parindent
        Department of Mathematics \\
        Oklahoma State University \\
        Stillwater, OK USA}
\email{segerman@math.okstate.edu}

%%% Not great -- how to get better positioning?
\thanks{This work is in the public domain.}

\date{\today}

%%% \begin{abstract}
%%% We discuss the art and science of producing conformally correct
%%% euclidean and hyperbolic tilings of compact surfaces.  As an example,
%%% we present a tiling of the Chmutov surface by hyperbolic (2, 4, 6)
%%% triangles.
%%% \end{abstract}

%%% MSC numbers: 51F15, 53A05, 49Q10, 65D17
%%% Keywords: conformal, Chmutov, tiling, meshes, uniformisation

\maketitle

\begin{figure}[htbp]
\includegraphics[width=0.72\textwidth]{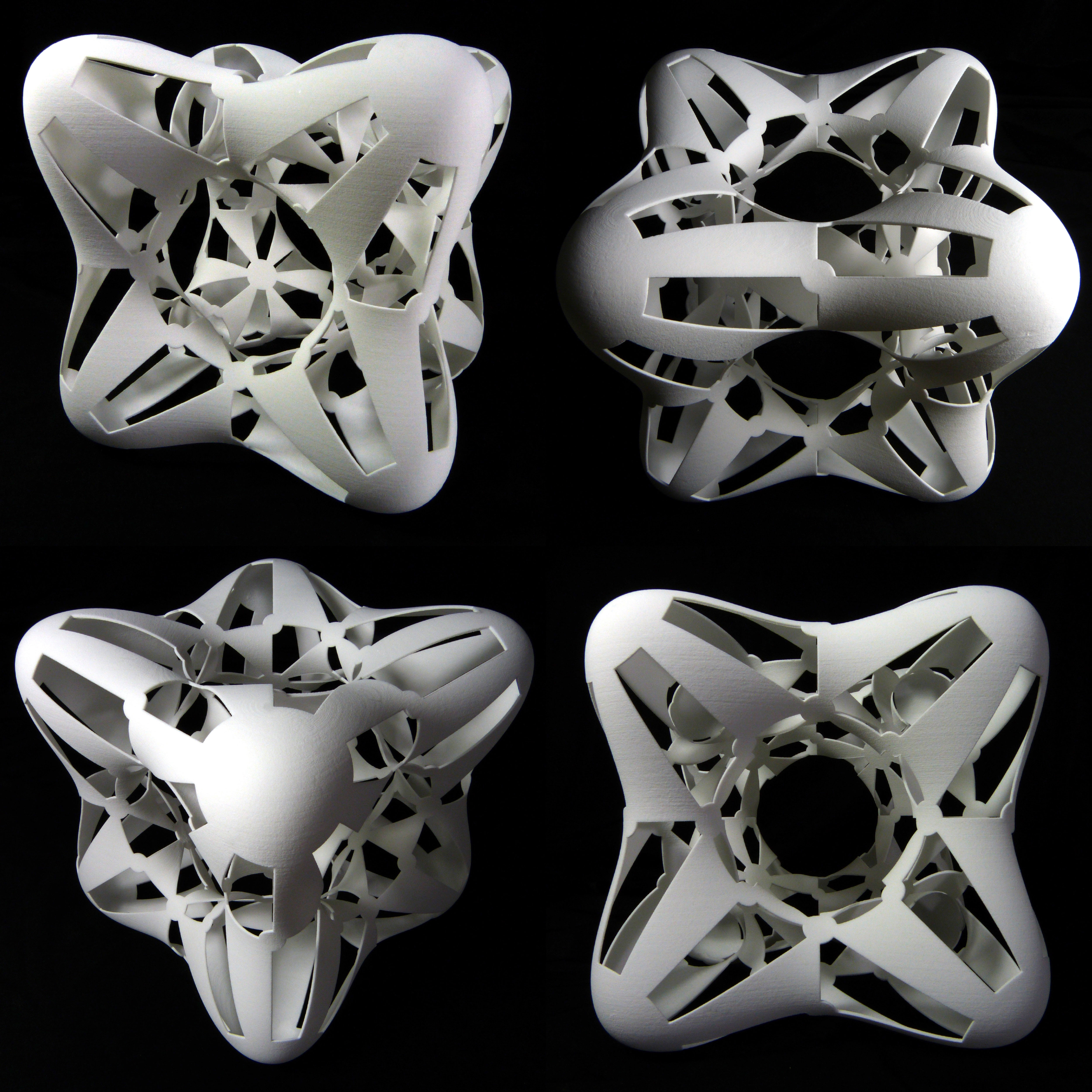}
\caption{A montage of photographs of a 3D print of the Conformal
  Chmutov sculpture. A 3D model can be viewed at
  SketchFab -- \url{https://skfb.ly/XYAC}.}

\label{Fig:Chmutov_training_montage}
\end{figure}

\section{Introduction}
\label{Sec:Intro}

Symmetric tilings are an ancient, ubiquitous, and beautiful motif in
decoration.  Perhaps the most famous examples are found at the
Alhambra in Granada.  Of course, a physical tiling on a wall or floor
ends at the corners of the room.  This is unfortunate; in our mind the
tiling goes on forever.  This raises the question of how we can
capture an infinite tiling in a finite space.

%%% O God, I could be bounded in a nut shell and count myself a king
%%% of infinite space, were it not that I have bad dreams.
%%%    -- Hamlet

When we try to wrap up something infinite, into something finite, we
often have to accept certain distortions of geometry.  The goal of
this note is to use modern tools to explore an old solution to this
problem: we give examples of \emph{conformally correct} tilings of
surfaces in three-space.

%%% ``Old solution'' refers to the uniformization theorem.

\section{Euclidean tilings}
\label{Sec:Euc}

\begin{wrapfigure}[14]{l}{0.4\textwidth}
  \vspace{-7pt}
  \centering
  \includegraphics[width=0.3\textwidth]{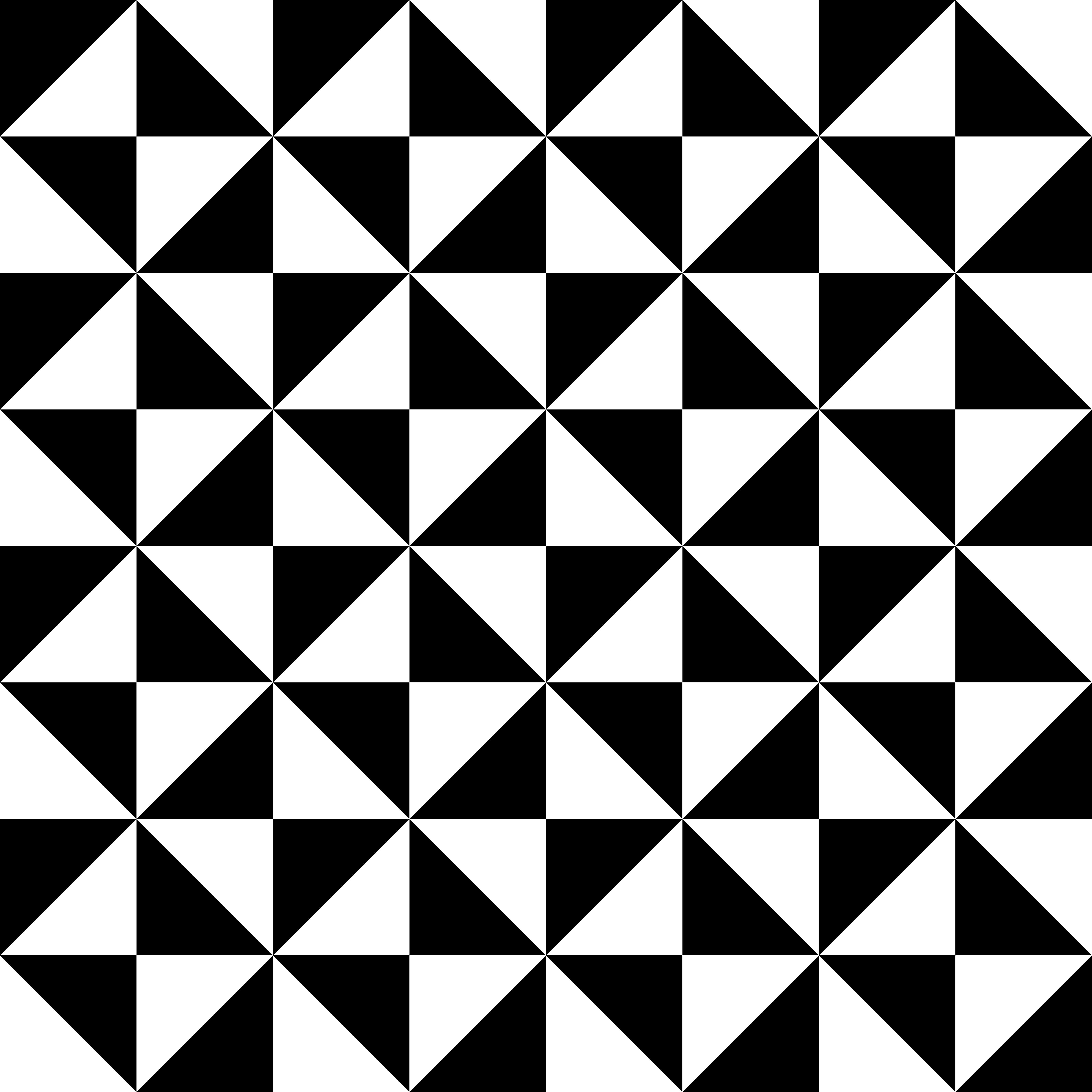}
  \caption{The (2,4,4) triangle tiling.} 
  \label{Fig:442_square}
\end{wrapfigure}

Perhaps the most familiar tiling is a brick tiling of a wall or
pathway.  Each brick meets six identical bricks, and the pattern
repeats in a regular way.  To add some flavour to our tilings, we will
insist on triangular tiles instead of rectangles.  Further, we will
require that the triangles alternate in colour between black and
white; each black triangle is the mirror image of any white
triangle. We fondly recall our teachers' words: ``The sum of the
angles of a triangle is $180$ degrees.''  This, and our desire for
mirror symmetry, implies that the possible angles for our triangle are
$(90, 45, 45)$, $(90, 60, 30)$, or $(60, 60, 60)$.  We call these the
$(2, 4, 4)$ triangle, the $(2, 3, 6)$ triangle, and the $(3, 3, 3)$
triangle respectively.  In \reffig{442_square} we have drawn a tiling
by copies of the $(2, 4, 4)$ triangle.  To explain these names, we
examine the vertices of the tiling (that is, any point where a pair of
lines cross).  There are three kinds of vertex: those where two black
triangles meet, those where four meet in a left windmill, and those
where four meet in a right windmill.

%%% Thus, the $(2, 4, 4)$ triangle is named after the number of copies
%%% about each vertex.

It may amuse you take a few minutes to draw a tiling by black and
white $(2, 3, 6)$ triangles, following the same rules.

\section{Wrapping up}
\label{Sec:Wrap}

\begin{wrapfigure}[21]{r}{0.4\textwidth}
\vspace{-23pt}
\centering
\subfloat[The $(2,4,4)$ tiling wrapped about a cylinder.]{
  \includegraphics[width=0.3\textwidth]{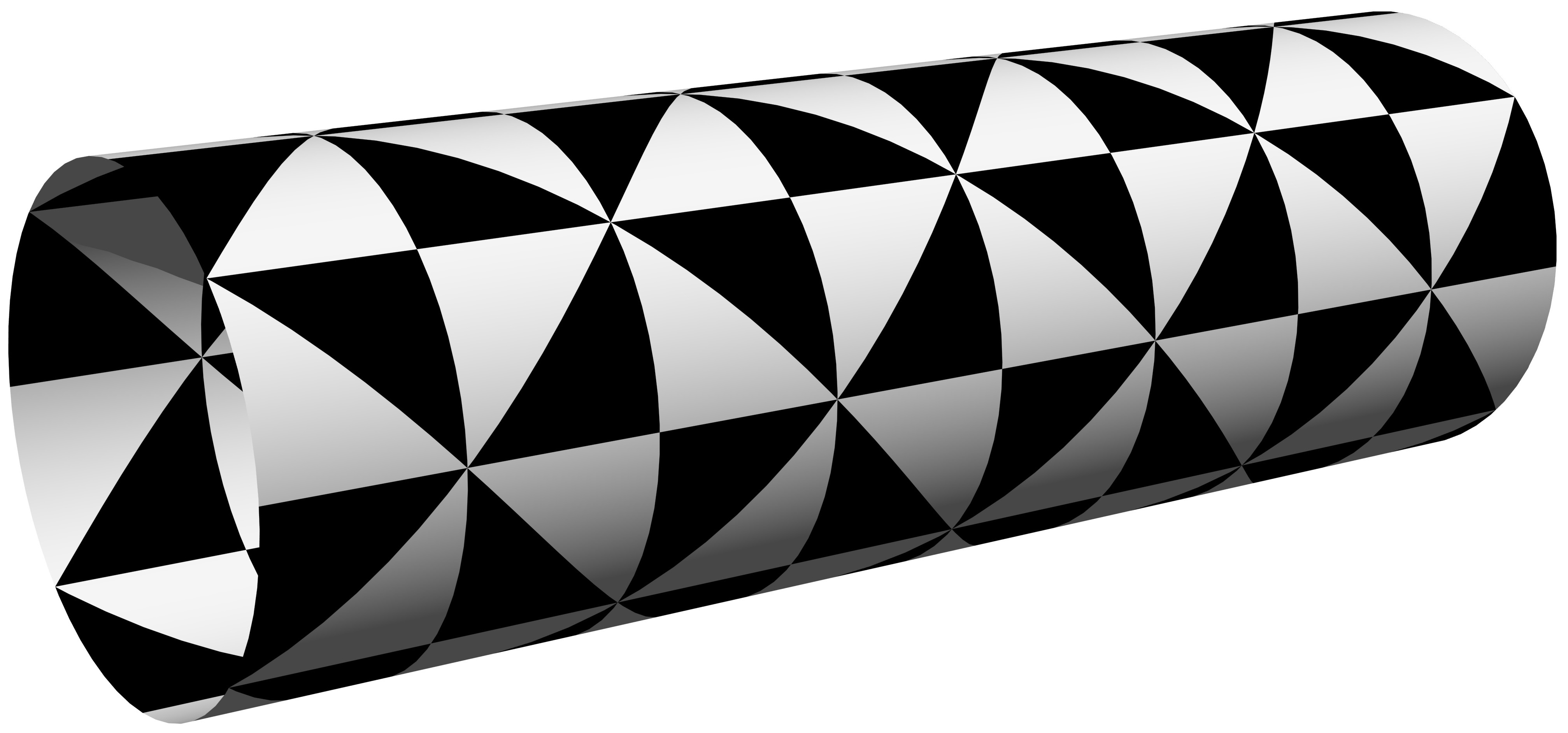}
  \label{Fig:442_cylinder}
}

\subfloat[The same tiling wrapped around a torus.]{
  \includegraphics[width=0.3\textwidth]{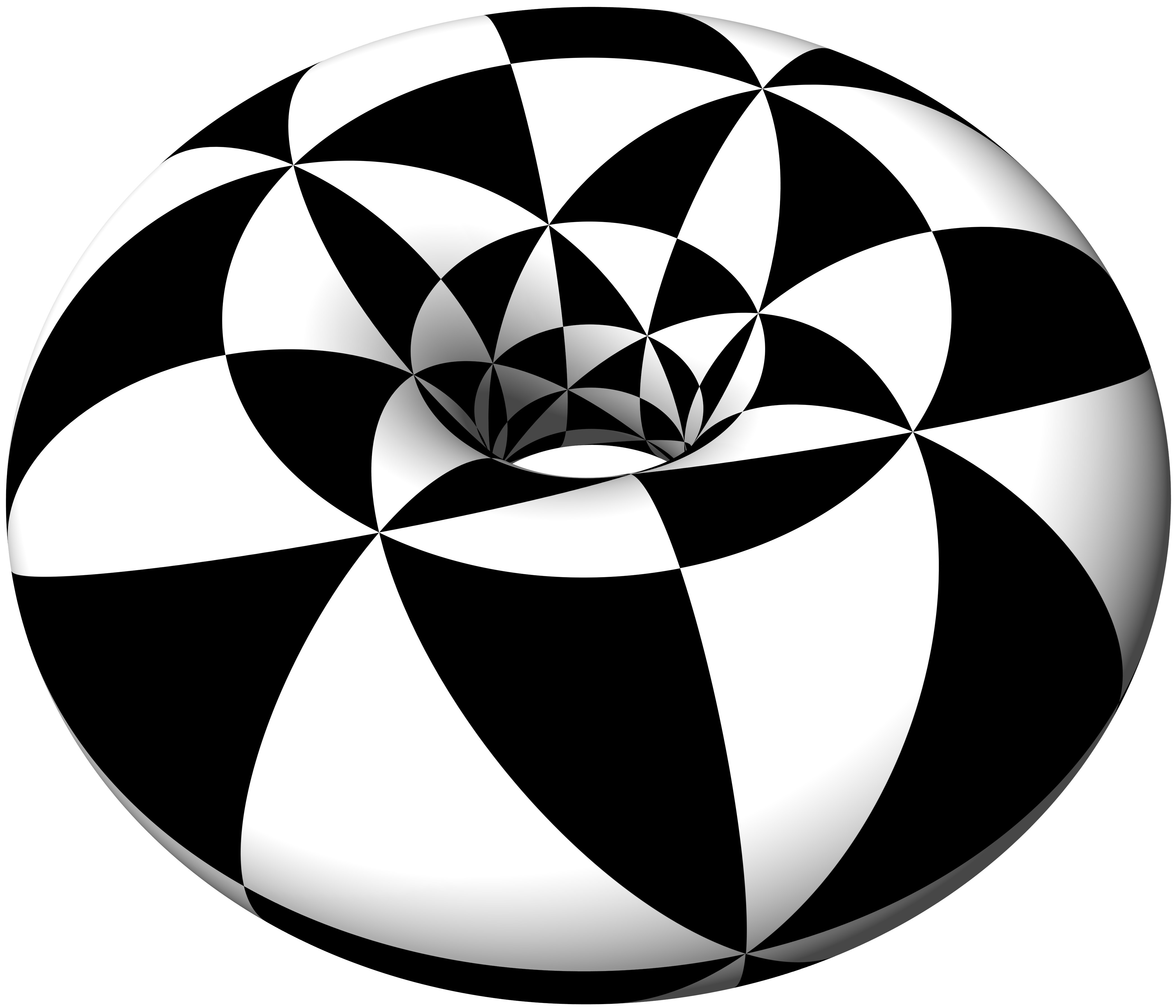}
  \label{Fig:442_torus}
}
\caption{Wrappings of the $(2, 4, 4)$ triangle tiling.} 
\label{Fig:Torus}
\end{wrapfigure}

In fact, the tiling in \reffig{442_square} is not infinite; it ends
rather abruptly, to avoid covering its caption.  This is an obvious
defect: it would be more pleasing if the tiling had no boundary edge.
Since the tiling cannot continue forever, we instead must find a way
to wrap it around something.

So, cut \reffig{442_square} out of the page, to get a tiled square.
Glue the top edge of the square to the bottom, forming a cylinder, as
shown in \reffig{442_cylinder}.  This gets rid of two of the four
edges of the square.  However, by bending the square to get a cylinder
we gave up something.  In the page all lines of the tiling are
straight; in the cylinder some lines become helices and some become
circles.  Nonetheless, the lengths of sides of triangles are unchanged
when we measure distance along the cylinder.  Likewise, the angles of
the triangles remain the same.

The next step of the wrapping process is much more difficult; we want
to glue the right boundary of the cylinder to the left.  It is an
interesting (but futile) exercise to try to do this second gluing with
a paper cylinder.  If we instead use some rubbery material, we get the
torus (the surface of a bagel) shown in \reffig{442_torus}.  This
solves our problem: the tiling is now finite but without boundary
edges.

%%% The need for the hole is perhaps a surprise, but it is required by
%%% Gauss-Bonnet; the average combinatorial curvature equals the
%%% average Gaussian curvature.

Again, to get something we had to give up something.  As they lie in
the torus (\reffig{442_torus}), the lengths of the sides of the
triangles are badly distorted: sides near the hole are squeezed while
those far from the hole are stretched.  This is unavoidable for the
following reason.  Near the hole the torus is negatively curved -- a
small part of it is saddle-shaped.  Far from the hole the torus is
positively curved -- a small part looks like a small part of a sphere.
And a $(2, 4, 4)$ tiling with correct lengths cannot lie snugly on a
surface where it has either positive or negative curvature.

However, if we are careful, some of the geometry of the tiling will
survive.  The vertical and horizontal lines of the tiling have become
circles that go around and go through the torus, respectively.  These
two families meet everywhere at right angles, just as in the original
tiling.  The diagonal lines, that were helices in the cylinder, have
become \emph{Villarceau circles}.  As first observed by
Sch{\oe}lcher~\cite{Mannheim03}, the Villarceau circles meet
everywhere with a constant angle.

If we stretch the torus to look like a bicycle tire then the
Villarceau angle decreases towards zero.  If we fatten the torus then
the angle increases to $180$ degrees.  So, as in the story of
Goldilocks and the Three Bears, there is a size, somewhere in between,
which is just right!
%%% IVT
At this special size, the Villarceau angle is $90$ degrees.
Since the torus has a reflection symmetry, at this special
size the Villarceau circles and the vertical circles meet at an angle
of $45$ degrees.  Thus the angles at the corners of the triangles are
exactly the same as in the tiling of the plane.
%%% Mannheim, M. A. (1903). ``Sur le th{\'e}or{\`e}me de
%%% Sch{\oe}lcher''. Nouvelles Annales de Mathématiques. Paris:
%%% Carilian-Gœury et Vor. Dalmont. 4th series, volume 3: 105–107.
%%%
%%% See also: 
%%%
%%% Yvon Villarceau, Antoine Joseph François (1848).
%%% ``Th{\'e}or{\`e}me sur le tore''.  Nouvelles Annales de
%%% Mathématiques. Série 1. Paris: Gauthier-Villars. 7: 345–347.
%%% OCLC: 2449182.
When this happens we say that the tiling is \emph{conformally correct}.

Now is a good time to eat a few oranges.  Once you are refreshed, you
might want to draw a conformally correct tiling by black and white
$(2, 3, 4)$ triangles.

%%% on one of your remaining oranges...

\section{Hyperbolic tilings}
\label{Sec:Hyp}

\begin{wrapfigure}[17]{l}{0.5\textwidth}
  \vspace{-15pt}
  \centering
  \includegraphics[width=0.4\textwidth]{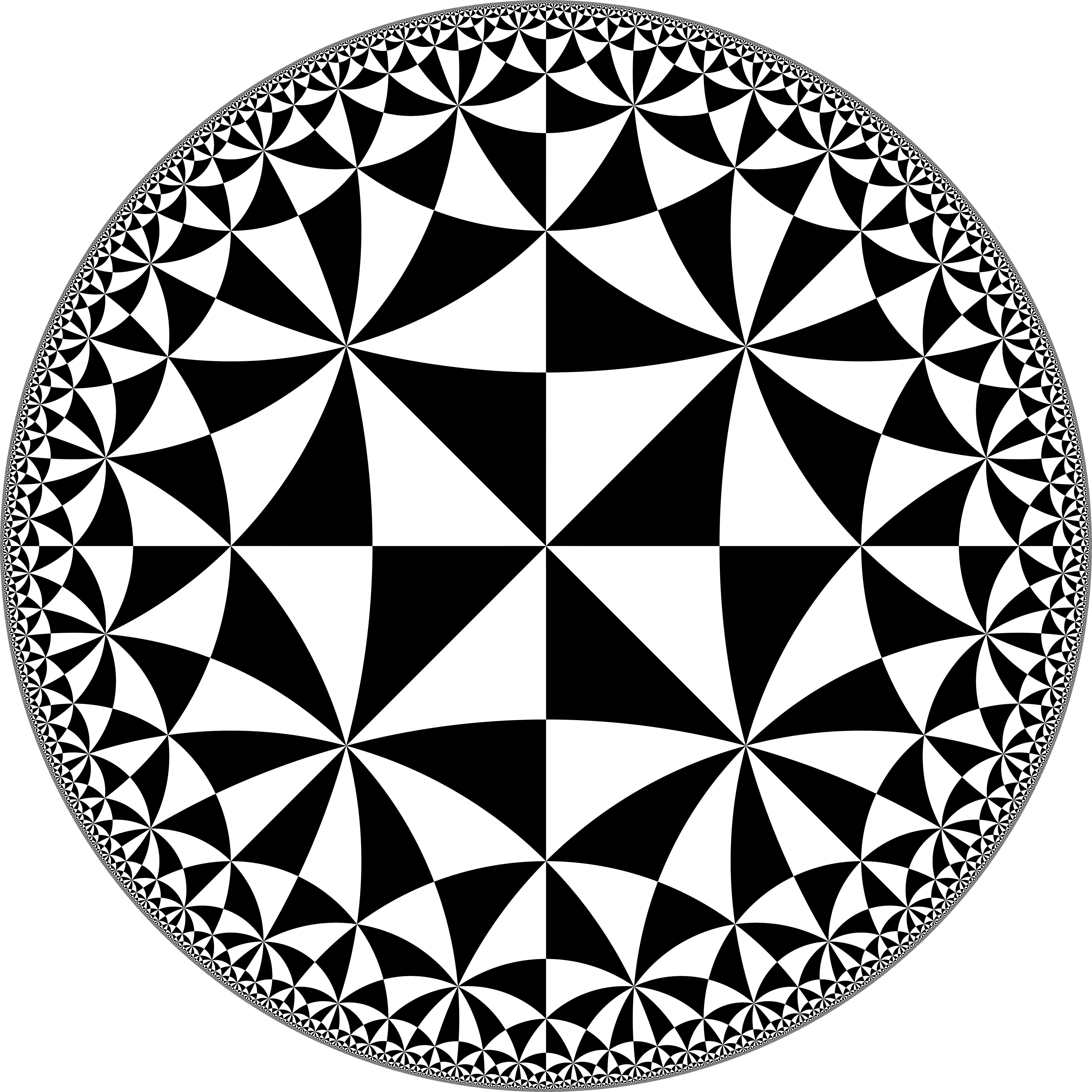}
  \caption{The $(2, 4, 6)$ triangle tiling of the hyperbolic plane,
    drawn by Roice Nelson.}
  \label{Fig:642}
\end{wrapfigure}

The three triangles $(2, 4, 4)$, $(2, 3, 6)$, and $(3, 3, 3)$ are the
only ones giving tilings of the flat, or \emph{euclidean}, plane.  If we
wish to find further examples we must learn to ignore our teachers; we
must consider the possibility of a triangle where the sum of the
angles is less than $180$ degrees.  These properly live in the
hyperbolic, or lobachevskian, plane.  In \reffig{642} we see a tiling
of the hyperbolic plane by identical (mirrored) black and white $(2,
4, 6)$ triangles.

You, gentle reader, have already noticed the obvious problem -- the
supposedly identical triangles are not identical at all.  The lengths
of edges shrink as the triangles march outwards.  This is an
unavoidable feature of hyperbolic geometry, first proved by
Hilbert~\cite{Hilbert01}: there is no picture of the hyperbolic plane
in euclidean three-space which faithfully represents lengths.
%%% rather, there is no smooth isometric immersion
Thus \reffig{642} is not a tiling of the hyperbolic plane but is
instead a \emph{model} of the tiling.

\begin{wrapfigure}[21]{r}{0.55\textwidth}
\vspace{-10pt}
\centering
\includegraphics[width=0.51\textwidth]{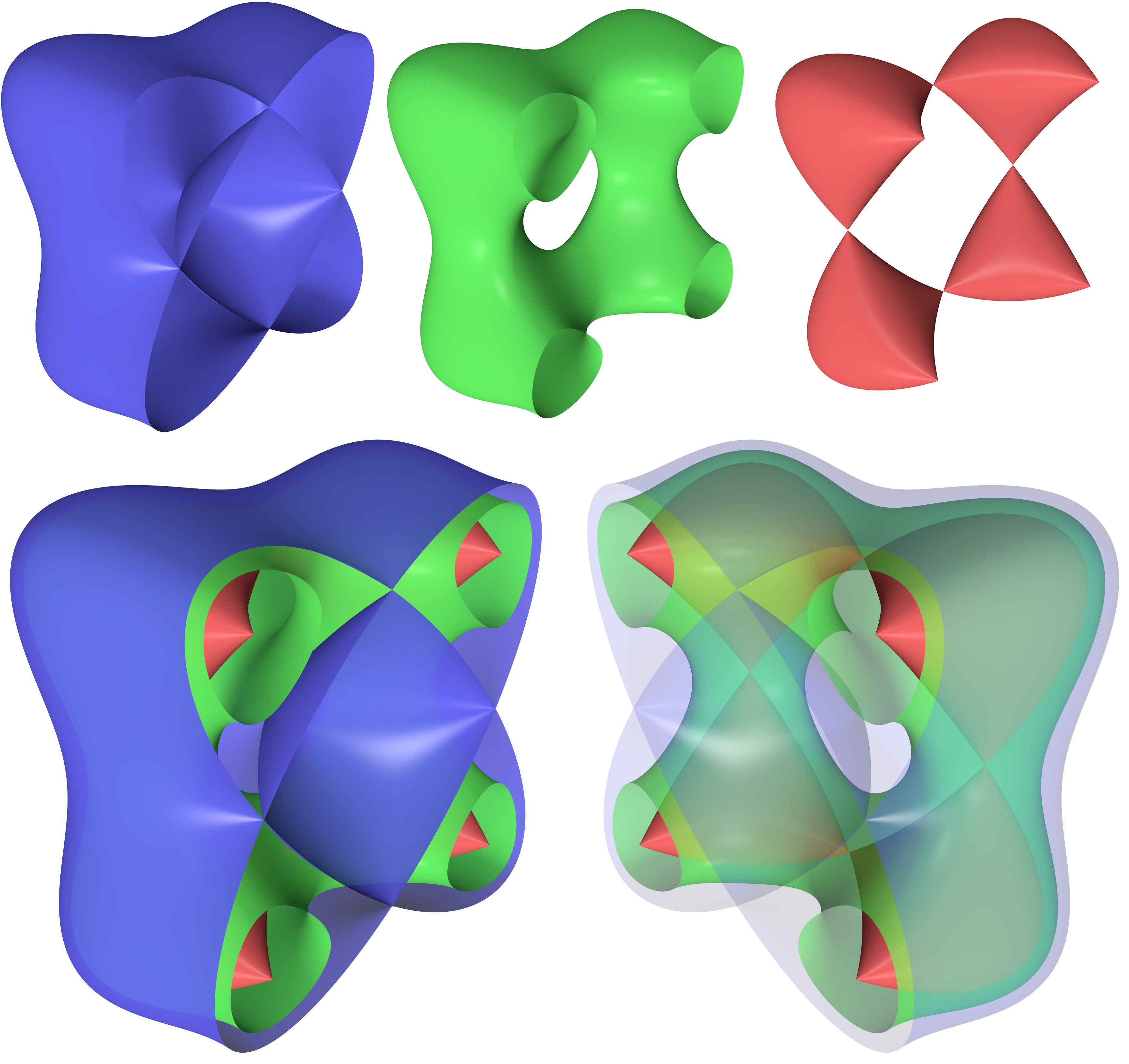}
\caption{Above: half of the surfaces $S(1)$ (in blue), $S(0)$ (in
  green), and $S(-1)$ (in red). Below left: the half surfaces arranged
  in space together, inspired by Curtis~\cite{Banchoff96}.  Below
  right: the same, drawn with greater transparency.}
\label{Fig:Chmutov_levels}
\end{wrapfigure}

We must content ourselves with the fact that this model is again
conformally correct: the angles are exact.  Consequently, the number
of black triangles about any vertex is either two, four, or six.

You may want to count the number of triangles at a fixed number of
triangle-sized steps away from the centre of the $(2, 4, 6)$ tiling.
This behaves very differently from the similar count in the $(2, 4,
4)$ tiling.

\section{Wrapping the hyperbolic plane}
\label{Sec:HypWrap}

%%% Email from Gary Kennedy: I've been working with Kate Shannon and
%%% John Thrasher in putting together the mathematical art exhibition
%%% that begins at Ohio State Mansfield on November 9th, although Kate
%%% has been doing most of the work, and in particular she's been
%%% communicating with our artists.  But I had an idea, and she
%%% suggested I write to you directly.  I'm wondering if we can
%%% custom-order a piece for our show, in addition to the ones you've
%%% already promised.  Specifically, I'm hoping we can persuade you to
%%% print a Chmutov surface.  Sergei Chmutov is a Professor of
%%% Mathematics at our campus.

In 2015, the Pearl Conard Art Gallery at the Ohio State University put
on an exhibition of mathematical art.  One of the organisers, Gary
Kennedy, asked us to contribute a sculpture based on the \emph{Chmutov
surface}, in honour of Professor Sergei Chmutov.  To understand this
surface, there is a necessary algebraic definition:
\[
F(x, y, z) = 8(x^4 + y^4 + z^4) - 8(x^2 + y^2 + z^2) + 3.
\]

For any real number $c$ we take $S(c)$ to be all of the points $(x, y,
z)$ in three-space solving the equation $F(x, y, z) = c$.  We call
these \emph{contour surfaces}.  The contour surfaces for $c = 1$, $c =
0$, and $c = -1$ are shown in \reffig{Chmutov_levels}.  At $c = 1$ the
contour has six nodes (sharp points at which the surface is not
smooth) corresponding to the faces of a cube; at $c = -1$ the contour
has twelve nodes corresponding to the edges of the same cube.  This
surface was first studied by Chmutov and Hirzebruch~\cite{Banchoff91}.
%%% The maximum possible number of nodes in degree four is 16 (the
%%% Kummer surface).

% \begin{wrapfigure}[39]{r}{0.35\textwidth}
\begin{figure}[htbp]
% \vspace{-25pt}
\centering
\subfloat[Half of the surfaces $S(c)$, with symmetry planes drawn in blue.]{
  \includegraphics[width=0.30\textwidth, valign=c]{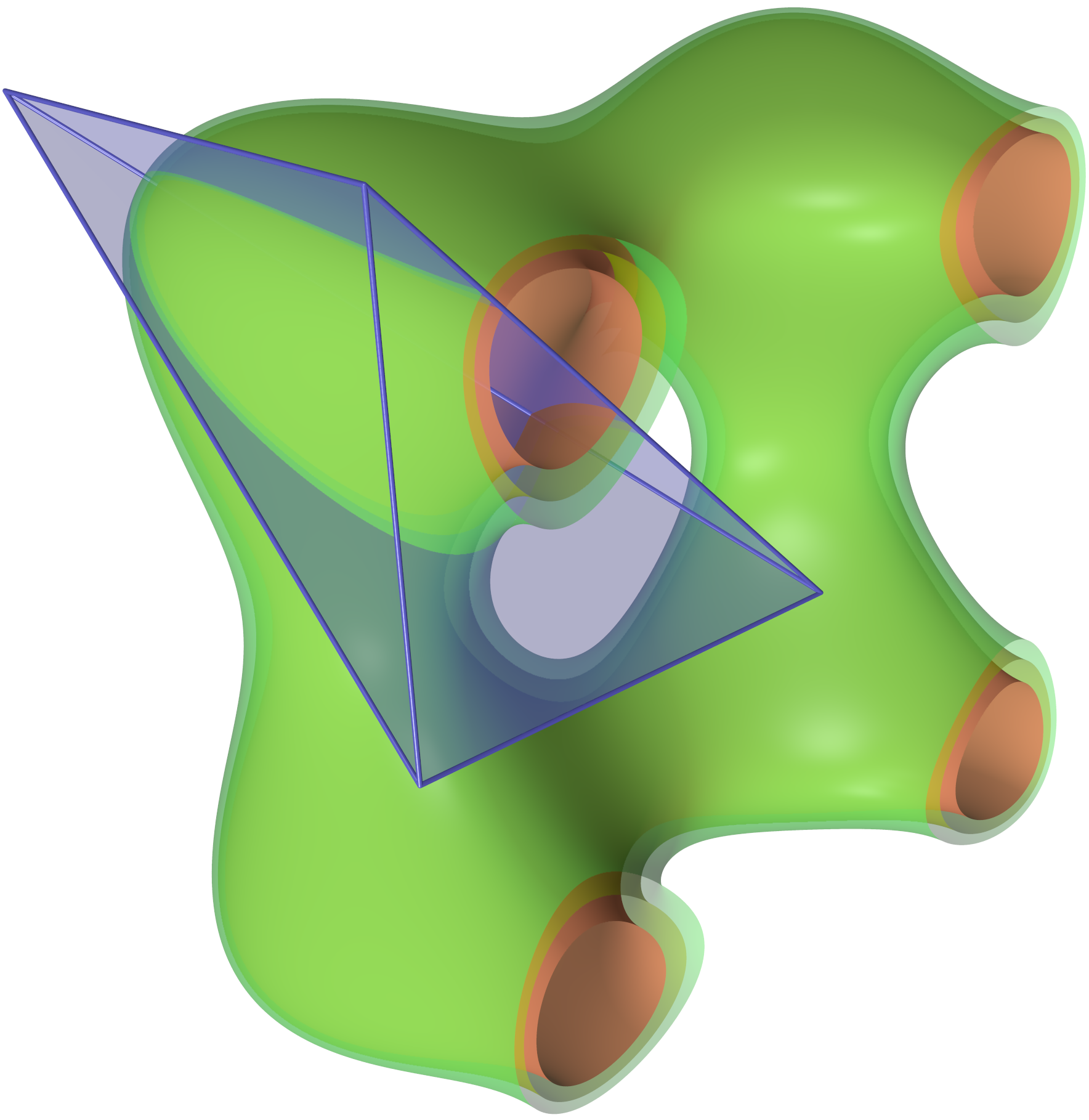}
  \label{Fig:chmutov_half_and_tet}
}
\subfloat[The patches $Q(c)$.]{
  \includegraphics[width=0.30\textwidth, valign=c]{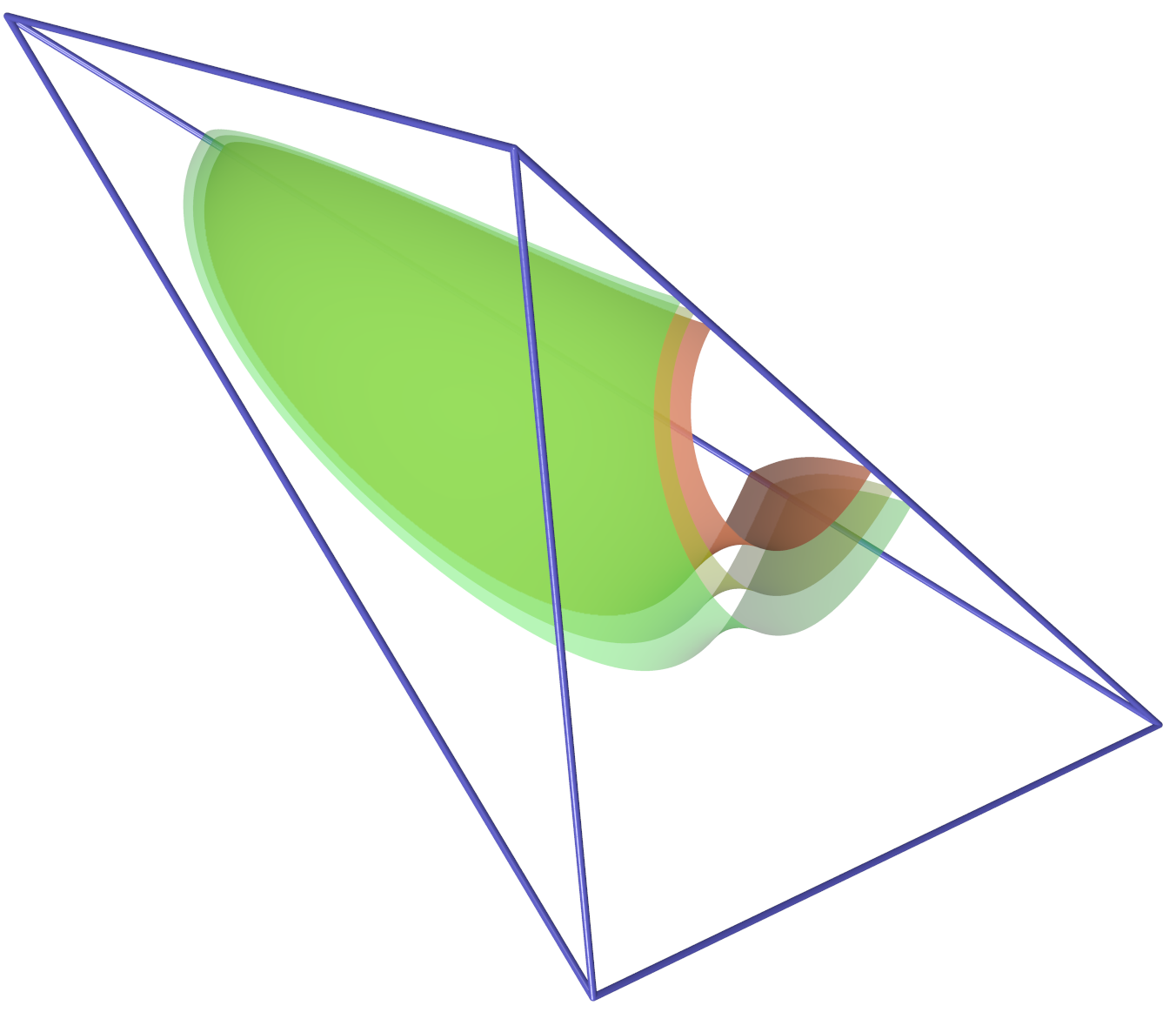}%
  \vphantom{\includegraphics[width=0.30\textwidth, valign=c]{chmutov_quads_R3_half_and_tet}}
  \label{Fig:chmutov_quads}
}
\subfloat[The patches laid out in the hyperbolic plane.]{
  \includegraphics[width=0.30\textwidth, valign=c]{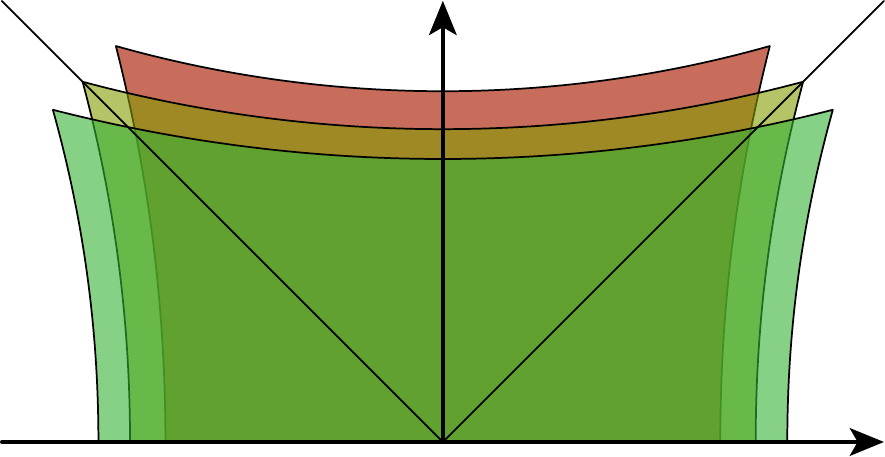}%
  \vphantom{\includegraphics[width=0.30\textwidth, valign=c]{chmutov_quads_R3_half_and_tet}}
  \label{Fig:quads_in_H2}
}
\caption{Parts of the surfaces $S(0)$ (in green), $S(-0.2411...)$ (in
  yellow), and $S(-0.5)$ (in orange).}
\label{Fig:search_chmutov}
\end{figure}

Our first idea for a sculpture was to wrap a hyperbolic tiling about
the contour $S(0)$; since the surface has more than one hole, no
spherical or euclidean tiling will work.  To decide which hyperbolic
tiling to use, first note that every contour $S(c)$ has the same
symmetries as the cube.  If we cut $S(c)$ along all symmetry planes of
the cube, as shown in \reffig{chmutov_half_and_tet}, we find a
\emph{patch} $Q(c)$ -- a four-sided polygon curving in space.  See
\reffig{chmutov_quads}.  By reflecting the patch around, we can
recover the whole surface.  At its four corners the patch $Q(c)$ has
angles $60$, $60$, $90$, and $90$ degrees.

According to the uniformisation theorem~\cite{StGervais16}
%%% really the Schwarz reflection principle
there is a unique way to lay $Q(c)$ out in the hyperbolic plane with
the correct angles at its corners.
%%% up to Mobius transformation
The program Confoo~\cite{vonGagern, SpringbornEtAl08} calculates how
to do this.  The boundary edges of the patches $Q(c)$ then receive
hyperbolic lengths.  If the vertical sides were exactly half the
length of the upper side, then we could divide $Q(c)$ into four $(2,
4, 6)$ triangles, two black and two white.  Unfortunately, for $Q(0)$
this is not the case!  In \reffig{chmutov_quads} we see that the green
rectangle is too wide.  If we increase $c$ from $0$, we get to the
blue surface in \reffig{Chmutov_levels}, which pinches off six nodes
at the faces of the cube.  For $Q(-0.5)$, the orange rectangle is too
tall. If we continue to decrease $c$ from $-0.5$, we get to the red
surface in \reffig{Chmutov_levels}, which pinches off twelve nodes at
the edges of the cube.  The surface gets too skinny near the faces in
one direction, and too skinny near the edges in the other.

So we can again apply the Goldilocks principle.  Using binary search
we find that $c \approx -0.2411$ gives a patch $Q(c)$ which is just
right.  The corresponding yellow rectangle shown in
\reffig{quads_in_H2} has exactly the same shape as the four black and
white triangles just above the centre of the $(2, 4, 6)$ tiling in
\reffig{642}.  After a significant further amount of work, we produced
the sculpture shown in \reffig{Chmutov_training_montage}.

It is amusing to count the number of triangles in the tiling.
%%% Euler characteristic.
%%% 1/2*(2 - 3 + 1/2 + 1/4 + 1/6)*192 = -8 = 2 - 2*5.
Is there a faster way than counting the triangles one at a time?

%%% In fact, every $(2, 4, 6)$ tiling of the Chmutov surface must have
%%% the same number of triangles.  This is very different from the
%%% case of the torus.

\subsection*{Acknowledgments} 

The Conformal Chmutov surface was commissioned by Gary Kennedy for the
Pearl Conard Art Gallery, at the Ohio State University.  We thank
Martin von Gagern and Boris Springborn for helping us with Confoo.  We
thank Roice Nelson for drawing \reffig{642}.  The other figures and 3D
prints were produced using Python, POVray, Rhino, and Shapeways.

\begin{figure}[htbp]
  \includegraphics[width=0.7\textwidth]{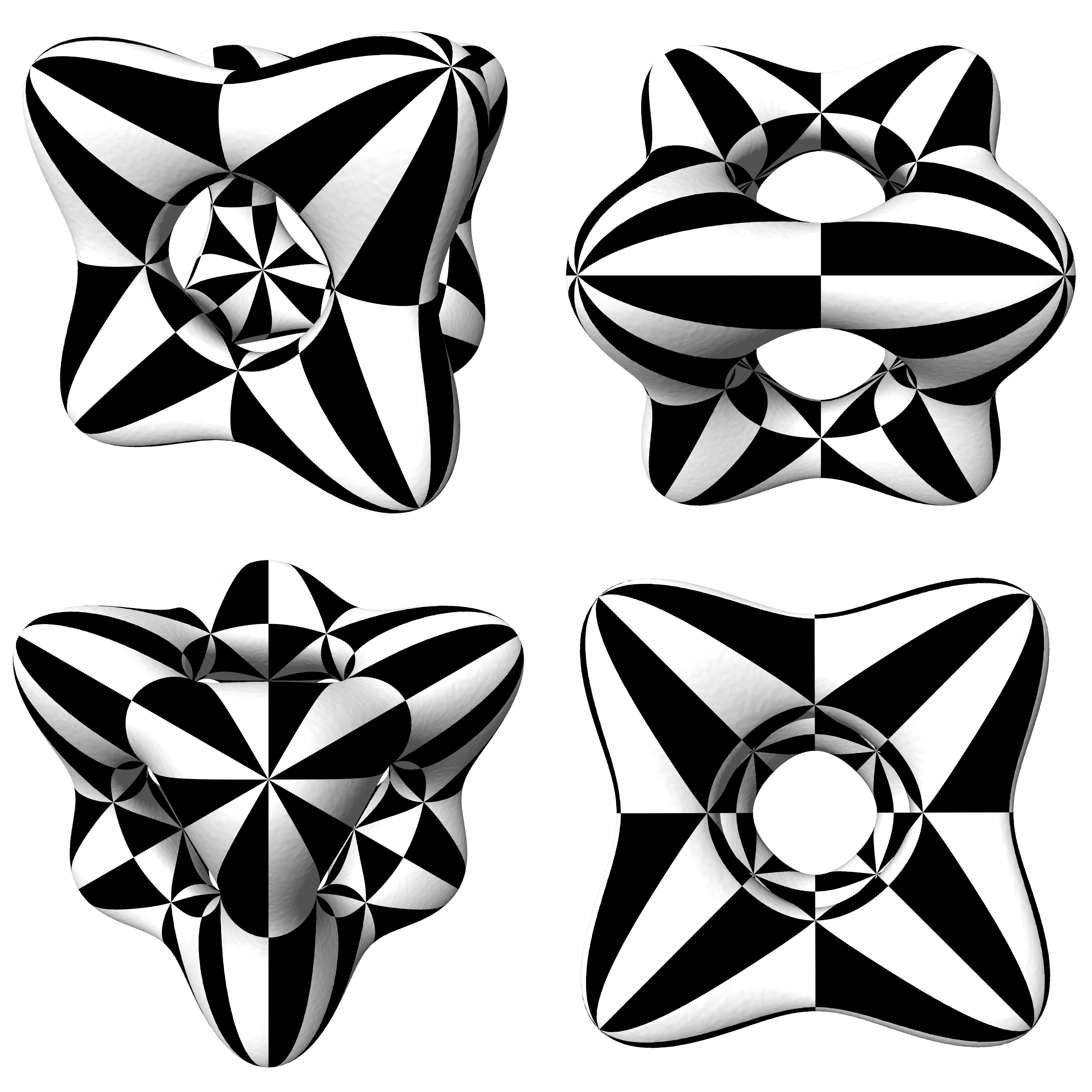}
  \caption{A montage of renders of a black and white tiling of the Conformal
    Chmutov sculpture.}
  \label{Fig:Chmutov_training_render_montage}
\end{figure}

\bibliographystyle{hyperplain}
\bibliography{bibfile}
\end{document}